\documentstyle{article}

\title{Coloured Hopf Algebras}
\author{C. Quesne\thanks{Directeur de recherches FNRS; E-mail:
cquesne@ulb.ac.be}\\ 
{\small Physique Nucl\'eaire Th\'eorique et Physique Math\'ematique,}\\ 
{\small Universit\'e Libre de Bruxelles, Campus de la Plaine CP229,}\\ 
{\small Boulevard~du Triomphe, B-1050 Brussels, Belgium}}
\date{}
\font\openface=msym10 

\def\R{\mbox{\openface R}}
\def\C{\mbox{\openface C}}
\def\Z{\mbox{\openface Z}}
\def\id{\mbox{\rm id}}
\def\case#1#2{{\textstyle{#1\over #2}}}

\newtheorem{definition}{Definition}[section]
\newtheorem{proposition}[definition]{Proposition}
\newtheorem{corollary}[definition]{Corollary}
\begin{document}

\maketitle

\section{Introduction}     
Since its introduction (for reviews see e.g.~\cite{majid}), the parametrized
Yang-Baxter equation (YBE) has played a crucial role in nonlinear integrable systems
in physics, such as exactly solvable statistical mechanics models and
low-dimensional integrable field theories. Its constant form has also
appeared in knot theory, where it is connected with braid groups. In addition,
the~YBE has inspired the development of quantum groups and quantum universal
enveloping algebras (QUEAs). The latter have appeared in the literature in two
different, but related forms: the Faddeev-Reshetikhin-Takhtajan (FRT) and 
Drinfeld-Jimbo (DJ) formulations.\par
%
%
In recent years, some integrable models with nonadditive-type solutions
$R^{\lambda,\mu} \ne R(\lambda - \mu)$ of the~YBE have been
discovered~\cite{bazhanov}. The corresponding~YBE
\begin{equation}
  R^{\lambda,\mu}_{12} R^{\lambda,\nu}_{13} R^{\mu,\nu}_{23} = R^{\mu,\nu}_{23}
  R^{\lambda,\nu}_{13} R^{\lambda,\mu}_{12}   \label{eq:colYBE}
\end{equation}
is referred to in the literature as the `coloured'~YBE, the nonadditive (in general
multicomponent) spectral parameters $\lambda$, $\mu$,~$\nu$ being considered
as `colour' indices. Constructing solutions of Eq.~(\ref{eq:colYBE}) has been
achieved by using various approaches (see
e.g.~\cite{burdik92a,gomez,kundu,bonatsos}). It should be stressed that this
coloured~YBE is distinct from the so-called `colour YBE'~\cite{mcanally} arising in
another context, as an extension of the graded~YBE to more general gradings than
that determined by $\Z_2$.\par
%
%
Extending the definitions of quantum groups and QUEAs by connecting
them to coloured $R$-matrices, instead of ordinary ones, has received some
attention in the literature. Kundu and Basu-Mallick~\cite{kundu} generalized the
FRT~formalism for some quantizations of $U(gl(2))$ and $Gl(2)$. In the context of
knot theory, Ohtsuki~\cite{ohtsuki} introduced some coloured quasitriangular Hopf
algebras, which are characterized by the existence of a coloured universal $\cal
R$-matrix, and he applied his theory to $U_q(sl(2))$. Bonatsos {\it et
al\/}~\cite{bonatsos} independently considered a rather similar, but nevertheless
distinct generalization for some nonlinear deformation of $U(su(2))$. Recently, we
extended the DJ formulation of QUEA's to coloured ones~\cite{cq96}, by elaborating
on the results of Bonatsos {\it et~al}.\par
%
%
It is the purpose of the present contribution to review such a generalization. In 
Sec.~2, we define coloured Hopf algebras in a way that generalizes Ohtsuki's
first attempt. In Sec.~3, we then apply the new concepts to construct two examples
of coloured QUEA's. Finally, Sec.~4 contains some concluding remarks.\par
%
%
\section{Coloured Hopf algebras}    
Let $\left({\cal H}_q, +, m_q, \iota_q, \Delta_q, \epsilon_q, S_q;k\right)$ (or in
short ${\cal H}_q$) be a Hopf algebra over some field~$k$ ($= \C$ or \R), depending
upon some parameters~$q$. Here $m_q: {\cal H}_q \otimes {\cal H}_q \to {\cal
H}_q$, $\iota_q: k \to {\cal H}_q$, $\Delta_q: {\cal H}_q \to {\cal H}_q \otimes
{\cal H}_q$, $\epsilon_q: {\cal H}_q \to k$, and $S_q: {\cal H}_q \to {\cal H}_q$
denote the multiplication, unit, comultiplication, counit, and antipode maps
respectively~\cite{majid}. Whenever $q$ runs over some set~$\cal Q$, called
{\it parameter set\/}, we obtain a set of Hopf algebras ${\cal H} = \{\,{\cal H}_q
\mid q \in {\cal Q}\,\}$. We may distinguish between two cases, according to
whether $\cal Q$ contains a single element (fixed-parameter case) or more than one
element (varying-parameter case).\par
%
%
Let us assume that there exists a set of one-to-one linear maps ${\cal G} = \{\,
\sigma^{\nu}: {\cal H}_q \to {\cal H}_{q^{\nu}} \mid q, q^{\nu} \in {\cal Q}, \nu \in
{\cal C}\,\}$, defined for any ${\cal H}_q \in {\cal H}$. They are labelled by some
parameters~$\nu$, called {\it colour parameters\/}, taking values in some
set~$\cal C$, called {\it colour set\/}. The latter may be finite, countably infinite,
or uncountably infinite. Two conditions are imposed on the $\sigma^{\nu}$'s:

\begin{itemize}
\item[(i)] Every $\sigma^{\nu}$ is an algebra isomorphism, i.e.,
\begin{equation}
  \sigma^{\nu} \circ m_q = m_{q^{\nu}} \circ \left(\sigma^{\nu} \otimes
  \sigma^{\nu}\right), \qquad \sigma^{\nu} \circ \iota_q = \iota_{q^{\nu}};
  \label{eq:iso} 
\end{equation}
\item[(ii)] $\cal G$ is a group (called {\it colour group\/}) with respect to the
composition of maps, i.e.,
\begin{eqnarray}
  \forall \nu, \nu' \in {\cal C}, \exists\, \nu'' \in {\cal C}: \sigma^{\nu''} & =
           &\sigma^{\nu'} \circ \sigma^{\nu}: {\cal H}_q \to {\cal H}_{q^{\nu''}} = 
          {\cal H}_{q^{\nu,\nu'}},  \label{eq:compo} \\
  \exists\, \nu^0 \in {\cal C}: \sigma^{\nu^0} & = & \id: {\cal H}_q \to 
          {\cal H}_{q^{\nu^0}} = {\cal H}_q, \\
  \forall \nu \in {\cal C}, \exists\, \nu' \in {\cal C}: \sigma^{\nu'} & =
           &\sigma_{\nu} \equiv \left(\sigma^{\nu}\right)^{-1}: {\cal H}_{q^{\nu}} \to
          {\cal H}_{q}   \label{eq:inverse}.
\end{eqnarray}
In Eqs.~(\ref{eq:compo}) and (\ref{eq:inverse}), $\nu''$ and $\nu'$ will be denoted by
$\nu' \circ \nu$ and $\nu^i$, respectively.
\end{itemize}
\par
%
%
$\cal H$, $\cal C$, and $\cal G$ can be combined into
%
\begin{definition}   \label{def-colmaps}
The maps $\Delta^{\lambda,\mu}_{q,\nu}: {\cal H}_{q^{\nu}} \to {\cal
H}_{q^{\lambda}} \otimes {\cal H}_{q^{\mu}}$, $\epsilon_{q,\nu}: {\cal H}_{q^{\nu}}
\to k$, and $S^{\mu}_{q,\nu}: {\cal H}_{q^{\nu}} \to {\cal H}_{q^{\mu}}$, defined by
\begin{equation}
  \Delta^{\lambda,\mu}_{q,\nu} \equiv \left(\sigma^{\lambda} \otimes
  \sigma^{\mu}\right) \circ \Delta_q \circ \sigma_{\nu}, \qquad
  \epsilon_{q,\nu} \equiv \epsilon_q \circ \sigma_{\nu}, \qquad
  S^{\mu}_{q,\nu} \equiv \sigma^{\mu} \circ S_q \circ \sigma_{\nu},
  \label{eq:colmaps} 
\end{equation}
for any $q \in {\cal Q}$, and any $\lambda$, $\mu$, $\nu \in {\cal C}$, are called
coloured comultiplication, counit, and antipode respectively. 
\end{definition}
\par
%
%
It is easy to prove the following proposition:
%
\begin{proposition}   \label{prop-GenHopf}
The coloured comultiplication, counit, and antipode maps, defined in
Eq.~(\ref{eq:colmaps}), transform under the colour group~$\cal G$ as
\begin{eqnarray}
  \left(\sigma^{\lambda}_{\alpha} \otimes \sigma^{\mu}_{\beta}\right) \circ 
           \Delta^{\alpha,\beta}_{q,\nu} & = & \Delta^{\lambda,\mu}_{q,\nu} = 
           \Delta^{\lambda,\mu}_{q,\gamma} \circ \sigma^{\gamma}_{\nu}, 
           \nonumber \\
  \epsilon_{q,\alpha} \circ \sigma^{\alpha}_{\nu} & = & \epsilon_{q,\nu},
           \nonumber \\
  \sigma^{\mu}_{\alpha} \circ S^{\alpha}_{q,\nu} & = & S^{\mu}_{q,\nu} = 
           S^{\mu}_{q,\beta} \circ \sigma^{\beta}_{\nu},
\end{eqnarray}
and satisfy generalized coassociativity, counit, and antipode axioms
\begin{eqnarray}
  \left(\Delta^{\alpha,\beta}_{q,\lambda} \otimes \sigma^{\gamma}_{\mu}\right)
           \circ \Delta^{\lambda,\mu}_{q,\nu} & = & \left(\sigma^{\alpha}_{\lambda'}
           \otimes \Delta^{\beta,\gamma}_{q,\mu'}\right) \circ 
           \Delta^{\lambda',\mu'}_{q,\nu}, \nonumber \\
  \left(\epsilon_{q,\lambda} \otimes \sigma^{\alpha}_{\mu}\right) \circ
           \Delta^{\lambda,\mu}_{q,\nu} & = & \left(\sigma^{\alpha}_{\lambda'}
           \otimes \epsilon_{q,\mu'}\right) \circ \Delta^{\lambda',\mu'}_{q,\nu} =
           \sigma^{\alpha}_{\nu}, \nonumber \\
  m_{q^{\alpha}} \circ \left(S^{\alpha}_{q,\lambda} \otimes \sigma^{\alpha}_{\mu}
           \right) \circ \Delta^{\lambda,\mu}_{q,\nu} & = & m_{q^{\alpha}} \circ \left(
           \sigma^{\alpha}_{\lambda'} \otimes S^{\alpha}_{q,\mu'} \right) \circ
           \Delta^{\lambda',\mu'}_{q,\nu} = \iota_{q^{\alpha}} \circ \epsilon_{q,\nu},
           \label{eq:colantipode}
\end{eqnarray}
as well as generalized bialgebra axioms
\begin{eqnarray}
  \Delta^{\lambda,\mu}_{q,\nu} \circ m_{q^{\nu}} & = & \left(m_{q^{\lambda}}
           \otimes m_{q^{\mu}}\right) \circ (\id \otimes \tau \otimes \id) \circ
           \left(\Delta^{\lambda,\mu}_{q,\nu} \otimes \Delta^{\lambda,\mu}_{q,\nu}
           \right), \nonumber \\
  \Delta^{\lambda,\mu}_{q,\nu} \circ \iota_{q^{\nu}} & = & \iota_{q^{\lambda}}
           \otimes \iota_{q^{\mu}}, \nonumber \\
  \epsilon_{q,\nu} \circ m_{q^{\nu}} & = & \epsilon_{q,\nu} \otimes \epsilon_{q,\nu},
           \nonumber \\
  \epsilon_{q,\nu} \circ \iota_{q^{\nu}} & = & 1_k.
\end{eqnarray}
Here $\sigma^{\lambda}_{\mu}$ is the element of~$\cal G$ defined by
\begin{equation}
  \sigma^{\lambda}_{\mu} \equiv \sigma^{\lambda} \circ \sigma_{\mu},
\end{equation}
$\tau$ is the twist map, i.e., $\tau(a \otimes b) = b \otimes a$, $1_k$ denotes the
unit of~$k$, and no summation is implied over repeated indices.
\end{proposition}
\par
%
%
{}From Proposition~\ref{prop-GenHopf}, it is straightforward to obtain
%
\begin{corollary}   \label{corol-newHopf}
If Eqs.~(\ref{eq:iso})--(\ref{eq:inverse}) are satisfied, then for any $q \in \cal Q$,
any $\nu \in \cal C$, and $q_{\nu} \equiv q^{\nu^i}$, $\left({\cal H}_q, +, m_q,
\iota_q, \Delta^{\nu,\nu}_{q_{\nu},\nu}, \epsilon_{q_{\nu},\nu},
S^{\nu}_{q_{\nu},\nu}; k\right)$ is a Hopf algebra over $k$ with comultiplication
$\Delta^{\nu,\nu}_{q_{\nu},\nu}$, counit $\epsilon_{q_{\nu},\nu}$, and antipode
$S^{\nu}_{q_{\nu},\nu}$, defined by particularizing Eq.~(\ref{eq:colmaps}).
\end{corollary}
%
{\it Remark.\/} In particular, for $\nu = \nu^0$, we get back the original Hopf
structure of~${\cal H}_q$.\par
%
%
Generalizing the result contained in Corollary~\ref{corol-newHopf}, we are led to
introduce
%
\begin{definition}
A set of Hopf algebras~$\cal H$, endowed with coloured comultiplication, counit,
and antipode maps $\Delta^{\lambda,\mu}_{q,\nu}$, $\epsilon_{q,\nu}$,
$S^{\mu}_{q,\nu}$, as defined in (\ref{eq:colmaps}), is called coloured Hopf algebra,
and denoted by any one of the symbols $\left({\cal H}_q, +, m_q, \iota_q,
\Delta^{\lambda,\mu}_{q,\nu}, \epsilon_{q,\nu}, S^{\mu}_{q,\nu}; k, {\cal Q}, {\cal
C}, {\cal G}\right)$, $\left({\cal H}, {\cal C}, {\cal G}\right)$, or ${\cal H}^c$.
\end{definition}
%
\par
%
%
Let us now assume that the members of the Hopf algebra set~$\cal H$ are almost
cocommutative Hopf algebras~\cite{majid}, i.e., for any $q \in \cal Q$ there exists
an invertible element ${\cal R}_q \in {\cal H}_q \otimes {\cal H}_q$ (completed
tensor product), such that
\begin{equation}
  \tau \circ \Delta_q(a) = {\cal R}_q \Delta_q(a) {\cal R}_q^{-1}
\end{equation}
for any $a \in {\cal H}_q$.\par
%
%
We may then introduce
%
\begin{definition}    \label{def-colR}
Let ${\cal R}^c$ denote the set of elements ${\cal R}^{\lambda,\mu}_q \in {\cal
H}_{q^{\lambda}} \otimes {\cal H}_{q^{\mu}}$, defined by
\begin{equation}
  {\cal R}^{\lambda,\mu}_q \equiv \left(\sigma^{\lambda} \otimes
  \sigma^{\mu}\right) \left({\cal R}_q\right),    \label{eq:colR}  
\end{equation}
where $q$ runs over~$\cal Q$, and $\lambda$, $\mu$ over~$\cal C$.
\end{definition}
%
The following result can be easily obtained:
%
\begin{proposition}
If the Hopf algebras~${\cal H}_q$ of~$\cal H$ are almost cocommutative, then
${\cal R}^{\lambda,\mu}_q$, as defined in (\ref{eq:colR}), is invertible with
$\left({\cal R}^{\lambda,\mu}_q\right)^{-1}$ given by
\begin{equation}
  \left({\cal R}^{\lambda,\mu}_q\right)^{-1} = \left(\sigma^{\lambda} \otimes
  \sigma^{\mu}\right) \left({\cal R}_q^{-1}\right),   \label{eq:alcocom1}
\end{equation}
and
\begin{equation}
  \tau \circ \Delta^{\mu,\lambda}_{q,\nu}(a) = {\cal R}^{\lambda,\mu}_q
  \Delta^{\lambda,\mu}_{q,\nu}(a) \left({\cal R}^{\lambda,\mu}_q\right)^{-1}
  \label{eq:alcocom2}
\end{equation}
for any $a \in {\cal H}_{q^{\nu}}$.
\end{proposition}
%
Hence we have
%
\begin{definition}
A coloured, almost cocommutative Hopf algebra is a pair $\left({\cal H}^c, {\cal
R}^c\right)$, where ${\cal H}^c$ is a coloured Hopf algebra, ${\cal R}^c = \{\, {\cal
R}^{\lambda,\mu}_q \mid q \in {\cal Q}, \lambda, \mu \in {\cal C}\,\}$, and ${\cal
R}^{\lambda,\mu}_q$, defined in (\ref{eq:colR}), satisfies Eqs.~(\ref{eq:alcocom1})
and~(\ref{eq:alcocom2}).
\end{definition}
\par
%
%
In the same way, we may define coloured, quasitriangular (or coboundary, or
triangular) Hopf algebras by starting with an appropriate~${\cal R}_q$:
%
\begin{definition}   \label{def-quasi}
A coloured, almost cocommutative Hopf algebra $\left({\cal H}^c, {\cal R}^c\right)$
is said to be quasitriangular if
\begin{eqnarray}
  \left(\Delta^{\alpha,\beta}_{q,\lambda} \otimes \sigma^{\gamma}_{\mu}\right)
            \left({\cal R}^{\lambda,\mu}_q\right) & = & {\cal R}^{\alpha,\gamma}_{q,13}
            {\cal R}^{\beta,\gamma}_{q,23}, \nonumber \\
  \left(\sigma^{\alpha}_{\lambda} \otimes \Delta^{\beta,\gamma}_{q,\mu}\right)
            \left({\cal R}^{\lambda,\mu}_q\right) & = & {\cal R}^{\alpha,\gamma}_{q,13}
            {\cal R}^{\alpha,\beta}_{q,12}.    \label{eq:quasi}
\end{eqnarray}
The set~${\cal R}^c$ is called the coloured universal $\cal R$-matrix of $\left({\cal
H}^c, {\cal R}^c\right)$.
\end{definition}
\par
%
%
The terminology used for~${\cal R}^c$ in Definition~\ref{def-quasi} is justified by
the following proposition, which shows among others that the
elements of~${\cal R}^c$ satisfy the coloured YBE, as given in~(\ref{eq:colYBE}):
%
\begin{proposition}
Let $\left({\cal H}^c, {\cal R}^c\right)$ be a coloured quasitriangular Hopf algebra.
Then
\begin{eqnarray}
  {\cal R}^{\lambda,\mu}_{q,12} {\cal R}^{\lambda,\nu}_{q,13} 
            {\cal R}^{\mu,\nu}_{q,23} & = & {\cal R}^{\mu,\nu}_{q,23}
            {\cal R}^{\lambda,\nu}_{q,13} {\cal R}^{\lambda,\mu}_{q,12}, \nonumber \\
  \left(\epsilon_{q,\lambda} \otimes \sigma^{\alpha}_{\mu}\right) \left({\cal
            R}^{\lambda,\mu}_q\right) & = & \left(\sigma^{\alpha}_{\lambda'} \otimes
            \epsilon_{q,\mu'}\right) \left({\cal R}^{\lambda',\mu'}_q\right) =
           1_{q^{\alpha}}, \nonumber \\
  \left(S^{\alpha}_{q,\lambda} \otimes \sigma^{\beta}_{\mu}\right) \left({\cal
            R}^{\lambda,\mu}_q\right) & = & \left(\sigma^{\alpha}_{\lambda'} \otimes
            \left(S^{\mu'}_{q,\beta}\right)^{-1}\right) \left({\cal
            R}^{\lambda',\mu'}_q\right) = \left({\cal
            R}^{\alpha,\beta}_q\right)^{-1},   \label{eq:propquasi}  
\end{eqnarray}
where $1_{q^{\alpha}}$ denotes the unit element of~${\cal H}_{q^{\alpha}}$, and
$\left(S^{\mu}_{q,\nu}\right)^{-1}: {\cal H}_{q^{\mu}} \to {\cal H}_{q^{\nu}}$ is
given by $\left(S^{\mu}_{q,\nu}\right)^{-1} = \sigma^{\nu} \circ S_q^{-1} \circ
\sigma_{\mu}$.
\end{proposition}
%
\par
%
%
\section{Examples of coloured quantum universal enveloping algebras} 
In the present section, we construct two examples of coloured
quasitriangular Hopf algebras, for which the underlying Hopf algebras~${\cal H}_q$
are QUEAs of Lie algebras $U_q(g)$.\par
%
%
\subsection{The two-parameter quantum algebra $U_{q,s}(gl(2))$} 
The first example deals with the two-parameter deformation of
$U(gl(2))$~\cite{schirrmacher}, whose universal $\cal R$-matrix was given in
Ref.~\cite{burdik92b}. Such an example is quite significant since $U_{q,s}(gl(2))$
and its corresponding quantum group have played an important role both in
generating some matrix solutions of the coloured YBE~\cite{burdik92a,kundu},
and in constructing a coloured extension of the FRT 
formalism~\cite{kundu}.\par
%
%
The quantum algebra $U_{q,s}(gl(2))$, for which $k = \C$ and $q$, $s \in \C
\! \setminus \! \{0\}$, is generated by four operators $J_3$, $J_{\pm}$, $Z$, with
commutation relations
\begin{equation}
  \left[J_3, J_{\pm}\right] = \pm J_{\pm}, \qquad \left[J_+, J_-\right] = 
  \left[2J_3\right]_q, \qquad \left[Z, J_3\right] = \left[Z, J_{\pm}\right] = 0,
  \label{eq:gl(2)-def}
\end{equation}
and coalgebra and antipode depending upon both parameters $q$ and~$s$.\par
%
%
Eq.~(\ref{eq:gl(2)-def}) is left invariant under the transformations
\begin{equation}
  \sigma^{\nu}\left(J_3\right) = J_3, \qquad  \sigma^{\nu}\left(J_{\pm}\right) =
  J_{\pm}, \qquad  \sigma^{\nu}\left(Z\right) = \nu Z,
\end{equation}
where $\nu \in {\cal C} = \C \! \setminus \! \{0\}$, and $(q^{\nu},s^{\nu}) = (q,s)$
(fixed-parameter case). Since $\nu' \circ \nu = \nu' \nu$, $\nu^0 = 1$, $\nu^i =
\nu^{-1}$, the colour group $\cal G$ is isomorphic to the abelian group  $Gl(1,\C)$.
The coloured maps and universal $\cal R$-matrix are easily obtained as
\begin{eqnarray}
  \Delta^{\lambda,\mu}_{q,s,\nu}\left(J_3\right) & = & J_3 \otimes 1 + 1 \otimes
           J_3, \quad \Delta^{\lambda,\mu}_{q,s,\nu}\left(Z\right) = \frac{\lambda}
           {\nu}\, Z \otimes 1 + \frac{\mu}{\nu}\, 1 \otimes Z, \nonumber \\
  \Delta^{\lambda,\mu}_{q,s,\nu}\left(J_{\pm}\right) & = & J_{\pm} \otimes q^{J_3}
           \left(\frac{s}{q}\right)^{\pm\mu Z} + q^{- J_3} (q s)^{\pm\lambda Z} \otimes
           J_{\pm}, \nonumber \\
  \epsilon_{q,s,\nu}(X) & = & 0, \quad X \in \{J_3, J_{\pm}, Z\}, \nonumber \\
  S^{\mu}_{q,s,\nu}\left(J_3\right) & = & - J_3, \quad
           S^{\mu}_{q,s,\nu}\left(Z\right) = - \frac{\mu}{\nu} Z, \quad
           S^{\mu}_{q,s,\nu}\left(J_{\pm}\right) = - q^{\pm1} s^{\mp 2\mu Z} J_{\pm},
           \nonumber \\
  {\cal R}^{\lambda,\mu}_{q,s} & = & q^{2 \left(J_3 \otimes J_3 - \lambda Z
           \otimes J_3 + \mu J_3 \otimes Z\right)} \sum_{n=0}^{\infty}
           \frac{\left(1 - q^{-2}\right)^n}{[n]_q!}\, q^{n(n-1)/2} \nonumber \\
  & & \mbox{} \times \left(q^{J_3} (q s)^{- \lambda Z} J_+\right)^n \otimes 
           \left(q^{-J_3} \left(\frac{s}{q}\right)^{\mu Z} J_-\right)^n.   
\end{eqnarray}
\par
%
%
The matrix representation of the coloured universal $\cal R$-matrix in any
finite-dimensional representation of $U_{q,s}(gl(2))$ provides us with a matrix
solution $R^{\lambda,\mu}_q$ of the coloured YBE~(\ref{eq:colYBE}). For instance, in
the two-dimensional representation of $U_{q,s}(gl(2))$,
\begin{eqnarray}
  D(J_3) & = & \case{1}{2} \left(\begin{array}{cc}
                      1 & 0 \\ 0 & -1
                      \end{array} \right), \quad 
         D(J_+) = \left(\begin{array}{cc}
                       0 & 1 \\ 0 & 0
                       \end{array} \right), \quad
         D(J_-) = \left(\begin{array}{cc}
                       0 & 0 \\ 1 & 0
                       \end{array} \right), \nonumber \\
  D(Z) & = & \left(\begin{array}{cc}
                  1 & 0 \\ 0 & 1
                  \end{array} \right), 
\end{eqnarray}
we get a (renormalized) $4 \times 4$ coloured $R$-matrix $R^{\lambda,\mu}_{q,s} 
\equiv q^{1/2} (D \otimes D) \left({\cal R}^{\lambda,\mu}_{q,s}\right)$, given by
\begin{equation}
  R^{\lambda,\mu}_{q,s} = \left(\begin{array}{cccc}
              q^{1-\lambda+\mu} & 0 & 0 & 0 \\[0.1cm]
              0 & q^{\lambda+\mu} & \left(q - q^{-1}\right) s^{-\lambda+\mu} & 0 
                   \\[0.1cm]
              0 & 0 & q^{-\lambda-\mu} & 0 \\[0.1cm]
              0 & 0 & 0 & q^{1+\lambda-\mu}
              \end{array} \right).   
\end{equation}
The latter coincides with the coloured $R$-matrix previously derived by
Burd\'\i k and Hellin\-ger~\cite{burdik92a} by considering $2 \times 2$
representations of $U_{q,s}(gl(2))$ characterized by different eigenvalues
$\lambda$,~$\mu$ of~$Z$.\par
%
%
\subsection{The standard quantum oscillator algebra $U^{(s)}_z(h(4))$}
The second example is the standard deformation $U^{(s)}_z(h(4))$ of the oscillator
algebra $U(h(4))$, which was first derived by contracting
$U_q(gl(2))$~\cite{celeghini}, then recently obtained in a more convenient
basis~\cite{ballesteros}. This quantum algebra has been used to construct
a solution of the coloured YBE connected with some link
invariants~\cite{gomez}.\par
%
%
$U^{(s)}_z(h(4))$ is generated by four operators $N$, $M$,~$A_{\pm}$, satisfying
the commutation relations
\begin{equation}
  \left[N, A_{\pm}\right] = \pm A_{\pm}, \qquad \left[A_-, A_+\right] =
  \frac{\sinh(zM)}{z}, \qquad \left[M, N\right] = \left[M, A_{\pm}\right] = 0.
  \label{eq:s-h(4)-def}
\end{equation}
Here we assume $k = \C$ and $z \in \C \! \setminus \! \{0\}$. The algebra
defining relations~(\ref{eq:s-h(4)-def}) are left invariant under the
transformations
\begin{equation}
  \sigma^{\nu}(N) = N, \quad \sigma^{\nu}(M) = \nu_+ \nu_- M, \quad
  \sigma^{\nu}(A_+) = \nu_+ A_+, \quad \sigma^{\nu}(A_-) = \nu_- A_-,
  \label{eq:s-h(4)-G}
\end{equation}
where $\nu \equiv (\nu_+, \nu_-)$, provided $z$ is changed into $z^{\nu} = \nu_+
\nu_- z$. Hence the parameter set is ${\cal Q} = \C \!\setminus\! \{0\}$
(varying-parameter case), the colour set is the cartesian product ${\cal C} = (\C
\!\setminus\! \{0\}) \times (\C \! \setminus \! \{0\})$, and the colour group is the
direct product group ${\cal G} = Gl(1,\C) \otimes Gl(1,\C)$.\par
%
%
The corresponding coloured maps and universal $\cal R$-matrix are given by
\begin{eqnarray}
  \Delta^{\lambda,\mu}_{z,\nu}(N) & = & N \otimes 1 + 1 \otimes N, \quad
          \Delta^{\lambda,\mu}_{z,\nu}(M) = \frac{\lambda_+\lambda_-}{\nu_+\nu_-}\,
          M \otimes 1 + \frac{\mu_+\mu_-}{\nu_+\nu_-}\, 1 \otimes M, \nonumber \\
  \Delta^{\lambda,\mu}_{z,\nu}(A_+) & = & \frac{\lambda_+}{\nu_+}\, A_+ \otimes 1 
          + \frac{\mu_+}{\nu_+}\, e^{-\lambda_+\lambda_- z M} \otimes A_+,
          \nonumber \\
  \Delta^{\lambda,\mu}_{z,\nu}(A_-) & = & \frac{\lambda_-}{\nu_-}\, A_- \otimes  
          e^{\mu_+\mu_- z M}+ \frac{\mu_-}{\nu_-}\, 1 \otimes A_-, \nonumber \\
  \epsilon_{z,\nu}(X) & = & 0, \quad X \in \{N, M, A_{\pm}\}, \nonumber \\
  S^{\mu}_{z,\nu}(N) & = & - N, \quad S^{\mu}_{z,\nu}(M) = -
          \frac{\mu_+\mu_-}{\nu_+\nu_-} M, \nonumber \\
  S^{\mu}_{z,\nu}(A_{\pm}) & = & - \frac{\mu_{\pm}}{\nu_{\pm}} A_{\pm}\,
          e^{\pm\mu_+\mu_- z M}, \nonumber \\
  {\cal R}^{\lambda,\mu}_z & = & \exp\{-\lambda_+ \lambda_- z M \otimes N\}
          \exp\{-\mu_+ \mu_- z N \otimes M\} \nonumber \\
  & & \mbox{} \times \exp\{2\lambda_- \mu_+ z A_- \otimes A_+\}. 
          \label{eq:s-h(4)-Hopf}
\end{eqnarray}
\par
%
%
In the $3 \times 3$ matrix representation of $U^{(s)}_z(h(4))$ defined by
\begin{eqnarray}
  D(N) & = & \left(\begin{array}{ccc}
        0 & 0 & 0 \\ 0 & 1 & 0 \\ 0 & 0 & 0
        \end{array}\right), \qquad
  D(M) = \left(\begin{array}{ccc}
        0 & 0 & 1 \\ 0 & 0 & 0 \\ 0 & 0 & 0
        \end{array}\right), \nonumber \\
  D(A_+) & = & \left(\begin{array}{ccc}
        0 & 0 & 0 \\ 0 & 0 & 1 \\ 0 & 0 & 0
        \end{array}\right), \qquad
  D(A_-) = \left(\begin{array}{ccc}
        0 & 1 & 0 \\ 0 & 0 & 0 \\ 0 & 0 & 0
        \end{array}\right),   \label{eq:s-h(4)-rep}
\end{eqnarray}
the coloured universal $\cal R$-matrix is represented by the $9 \times 9$ matrix
\begin{equation}
  R^{\lambda,\mu}_z \equiv (D \otimes D)\left({\cal R}^{\lambda,\mu}_z\right) =
  \left(\begin{array}{ccc}
          1_3 & 2\lambda_-\mu_+ z D(A_+) & -\lambda_+\lambda_- z D(N) \\[0.1cm]
          0_3 & 1_3 - \mu_+\mu_- z D(M) & 0_3 \\[0.1cm]
          0_3 & 0_3 & 1_3
  \end{array}\right),   \label{s-h(4)-R}
\end{equation}
where $1_3$ and $0_3$ denote the $3\times 3$ unit and null matrices
respectively.\par
%
%
\section{Conclusion}
In this contribution, we did present some new algebraic structures, termed
coloured Hopf algebras~\cite{cq96}, by combining the coalgebra structures and
antipodes of a standard Hopf algebra set with the transformations of an algebra
isomorphism group, called colour group. We did also show that quasitriangular
Hopf algebras can be extended into coloured ones, characterized by the
existence of a coloured universal $\cal R$-matrix, satisfying the coloured YBE.
Finally, we did apply the new concepts to QUEAs of Lie algebras by constructing
two examples of coloured QUEAs.\par
%
%
As shown elsewhere~\cite{cq96}, the coloured Hopf algebras defined here
generalize those previously introduced by Ohtsuki~\cite{ohtsuki}, which are
restricted to abelian colour groups, in which case they reduce to substructures of
the present ones. Many more examples of coloured QUEAs, corresponding to
finite or infinite, abelian or nonabelian colour groups have been
constructed~\cite{cq96}, thereby providing some new solutions of the coloured
YBE, which might be of interest in the context of integrable models.\par
%
%
Extending the present formalism to QUEAs of Lie superalgebras, as well as defining
and constructing duals to coloured Hopf algebras are under current
investigation~\cite{cq97}.\par
%
%

\end{document}